\documentclass[leqno]{article}

\usepackage{amssymb}
\usepackage{amsfonts}
\usepackage{latexsym}
\usepackage[T1,T2A]{fontenc}
\usepackage[cp1251]{inputenc}
\usepackage[T2B]{fontenc}
\usepackage[bulgarian,russian,english]{babel}

\newtheorem{theorem}{Theorem}

\title {On some numerical characteristics of a bipartite graph
 }
\author {Krasimir Yordzhev}
\date {}

\begin {document}
\inputencoding{cp1251}

\maketitle
\begin{abstract}
The paper consider an equivalence relation in the set of vertices of a bipartite graph. Some numerical characteristics showing the cardinality of equivalence classes  are introduced. A combinatorial identity that is in relationship to these characteristics of the set of all  bipartite graphs of the type $g=\langle R_g \cup C_g , E_g \rangle$ is formulated and proved, where  $V=R_g \cup C_g$ is the set of vertices, $E_g$ is the set of edges of the graph $g$, $ |R_g |=m\ge 1$, $|C_g |= n\ge 1$, $|E_g |=k\ge 0$, $m,n$ and $k$ are integers.
\end{abstract}

\section{Introduction}
It is well known widespread use of graph theory in different areas of science and technology. For example, graph theory is a good tools for the modelling of computing devices and computational processes. So a lot of graph algorithms have been developed \cite{Mirchev,Swami}. One of the latest applications of the graph theory is to calculate the number of all disjoint pair of S-permutation matrices \cite{YordzhevISRN,Yordzhev2013}.
The concept of disjoint S-permutation matrices was introduced by Geir Dahl  \cite{dahl} in relation to the popular Sudoku puzzle.
 On the other hand, Sudoku matrices are special cases of Latin squares in the class of gerechte designs \cite{Bailey}.

Let $p$ be a positive integer. By $[p]$ we denote the set
$$[p] =\left\{ 1,2,\ldots ,p\right\} .$$

\emph{Bipartite graph} is the ordered triplet
 $$
g=\langle R_g \cup C_g , E_g \rangle ,
$$
where $R_g$ and $C_g$ are sets such that $R_g \ne \emptyset$, $C_g \ne \emptyset$ and $R_g \cap C_g =\emptyset$. The elements of the set
$$V_g =R_g \cup C_g$$
will be called \emph{vertices}. The set
$$E_g \subseteq R_g \times C_g =\{ \langle r,c \rangle \; |\; r\in R_g ,c\in C_g \}$$
will be called the set of \emph{edges}. Repeated edges are not allowed in our considerations.

Let $g'=\langle R_{g'} \cup C_{g'} , E_{g'} \rangle$ and $g''=\langle R_{g''} \cup C_{g''} , E_{g''} \rangle$. We will say that the graphs $g'$ and $g''$ are \emph{isomorphic} and we will write $g'\cong g''$, if $R_{g'} =R_{g''}$, $C_{g'} =C_{g''}$,  $|R_{g'}|=|R_{g''}|=m$, $|C_{g'}|=|C_{g''}|=n$ and there are $\rho \in \mathcal{S}_m$ and $\sigma \in\mathcal{S}_n$, where $\mathcal{S}_p$ is \emph{the symmetric group}, such that $\langle r, c\rangle \in E_{g'} \Longleftrightarrow \langle \rho (r) ,\sigma (c)\rangle \in E_{g''}$. The object of this work is bipartite graphs considered to within isomorphism.

Let $m,n$ and $k$ be integers, $m\ge 4$, $n\ge 1$, and let $0\le k\le mn$. Let us denote by $\mathfrak{G}_{m,n,k} $ the set of all  bipartite graphs without repeated edges of the type $g=\langle  R_g \cup C_g ,E_g \rangle $, considered to within isomorphism, such that $ |R_g |=m$, $|C_g |= n$ and $|E_g |=k$.

For more details on graph theory see  \cite{diestel,harary,Mirchev}.

In \cite{Fontana} Roberto Fontana proposed an algorithm which randomly gets a family of
$n^2 \times n^2$ mutually disjoint S-permutation matrices, where $n = 2, 3$. In $n = 3$ he
ran the algorithm 1000 times and found 105 different families of nine mutually
disjoint S-permutation matrices. Then, he obtained $9! \cdot 105 = 38\; 102\; 400$ Sudoku matrices.
In relation with Fontana's algorithm, it looks useful to calculate
the probability of two randomly generated S-permutation matrices to be disjoint.

The solution of this problem is given in \cite{Yordzhev2013}, where is described a formula for calculating all pairs of mutually disjoint S-permutation matrices. The application of this formula when $n=2$ and $n=3$ is explained in detail in  \cite{YordzhevISRN}.

To do that, the graph theory techniques have been used. It has been shown that to count the number of disjoint pairs of $n^2 \times n^2$ S-permutation matrices, it is sufficient to obtain some numerical characteristics of the set $\mathfrak{G}_{n,n,k}$ of all  bipartite graphs of the type $g=\langle R_g \cup C_g , E_g \rangle$, where  $V_g =R_g \cup C_g$ is the set of vertices, and $E_g$ is the set of edges of the graph $g$, $R_g \cap C_g =\emptyset$, $|R_g |=|C_g |=n$, $|E_g |=k$.

The aim of this work is to formulate and to prove a combinatorial problem, that is in relationship to some numerical characteristics of the elements of the set $\mathfrak{G}_{n,n,k}$.

For the classification  of all non defined concepts and notations
as well as for common assertions which have not been proved here see \cite{aigner,diestel,Sachkov}.

\section{An equivalence relation in a bipartite graph}

Let $$g=\langle R_g , C_g ,E_g \rangle \in \mathfrak{G}_{m,n,k}$$ for some natural numbers $m,n$ and $k$ and let $v\in V_g =R_g \cup C_g$.

By $N (v)$ we denote the set of all vertices of $V_g$, adjacent with $v$, i.e., $u\in N (v)$ if and only if there is an edge in $E_g$ connecting $u$ and $v$. In other words if $v\in R_g$, then $N(v) =\{ u\in C_g \; |\; \langle v,u \rangle \in E_g \}$ and if $v\in C_g $, then $N(v) =\{ u\in R_g \; |\; \langle u,v \rangle \in E_g \}$.
If $v$ is an isolated vertex  (i.e., there is no edge, incident with $v$), then by definition $N (v)=\emptyset$ and $ \textrm{degree}  (v) = |N(v)|=0$.

Since in $g$ there are not repeated edges, then it is easy to see that
$$\sum_{u\in R_{g}} |N (u)|=k \quad \& \quad \sum_{v\in C_{g}} |N (v)|=k \quad \Longrightarrow \quad \sum_{w\in V_{g}} |N (w)|=2k.$$

Let $g=\langle R_g , C_g ,E_g \rangle \in \mathfrak{G}_{m,n,k}$ and let $u,v\in V_g =R_g\cup C_g$. We will say that $u$ and $v$ are equivalent and we will write $u\sim v$ if $N (u) =N (v)$. If $u$ and $v$ are isolated, then by definition $u\sim v$ if and only if $u,v\in R_g$ or $u,v\in C_g$. Obviously if $u\sim v$, then  $u\in R_g \Leftrightarrow v\in R_g$ and $u\in C_g \Leftrightarrow v\in C_g$. It is easy to see that the above introduced relation is an equivalence relation.

By ${V_g}_{/\sim} $ we denote the obtained factor-set (the set of the equivalence classes) according to relation $\sim$ and let
$${V_g}_{/\sim} =\left\{ \Delta_1 ,\Delta_2 ,\ldots ,\Delta_s \right\} ,$$
where $\Delta_i \subseteq R_g$, or $\Delta_i \subseteq C_g$, $i=1,2,\ldots s$, $2\le s \le 2n$. We assume
$$\delta_i =|\Delta_i |,\quad 1\le \delta_i \le n , \quad i=1,2,\ldots , s$$
and for every $g\in \mathfrak{G}_{m,n,k}$ we define multi-set (set with repetition)
$$\left[ g \right] =\left\{ \delta_1 ,\delta_2 ,\ldots \delta_s \right\} ,$$
where $\delta_1 ,\delta_2 ,\ldots ,\delta_s$ are natural numbers, obtained by the above described way.

Obviously
$$\sum_{i=1}^s \delta_i =m+n .$$

The next assertion is a generalization of Corollary 1 of Lemma 1 from \cite{Yordzhev2013}.

\begin{theorem}\label{lrl1}
For every positive integers $m,n$ and every nonnegative integer $k$ such that $0\le k\le mn$ the following equation is true:
$$
\sum_{g\in \mathfrak{G}_{m,n,k} } \frac{1}{\displaystyle \prod_{\delta \in [g]} \delta !} =\frac{(mn)!}{m!n!k!(mn -k)!}
$$
\end{theorem}

Proof. A binary (or   boolean, or (0,1)-matrix)  is a matrix all of whose elements belong to the set $\mathfrak{B}
=\{ 0,1 \}$. With $b(m,n,k)$ we will denote the number of all  $m \times n$  binary matrices with exactly $k$ in number 1's, $k=0,1,\ldots ,mn$.

It is easy to see that
\begin{equation}\label{mnbink1}
b(m,n,k)={mn\choose k} = \frac{(mn)!}{k!(mn-k)!}
\end{equation}

We will prove that
\begin{equation}\label{bnk}
b(m,n,k)=m!n! \sum_{g\in \mathfrak{G}_{m,n,k} } \frac{1}{\displaystyle \prod_{\delta \in [g]} \delta !}
\end{equation}

Let $A=\left[ a_{ij} \right]_{m\times n}$ be $m\times n$ binary matrix with exactly $k$ 1's. Then we construct graph $g=\langle R_g\cup C_g ,E_g\rangle $, such that the set $R_g=\{ r_1 ,r_2 ,\ldots ,r_m \}$  corresponds to the rows of $A$, and $C_g=\{ c_1 , c_2 ,\ldots , c_n \}$  corresponds to the columns of $A$, however there is an edge connecting the vertices $r_i$ and $c_j$ if and only if $a_{ij} =1$. The graph, which has been constructed, obviously belongs to $\mathfrak{G}_{m,n,k}$.

Conversely, let $g=\langle R_g\cup C_g ,E_g\rangle \in \mathfrak{G}_{m,n,k}$. We number at a random way the vertices of $R_g$ by natural numbers from 1 to $m$ without repeating any of the numbers. This can be made by $m!$ ways. We analogously number the vertices of $C_g$  by natural numbers from 1 to $n$. This can be made by $n!$ ways. Then we construct the binary $m\times n$ matrix $A=\left[a_{ij} \right]_{m\times n}$, such that $a_{ij} =1$ if and only if there is an edge in $E_g$ connecting the vertex with number $i$ of $R_g$ with the vertex with number $j$ of $C_g$. As $g\in \mathfrak{G}_{m,n,k}$, then the matrix, that has been constructed, has exactly $k$ 1's. It is easy to see that when $q,r\in [m]$, $q$-th  and $r$-th rows of $A$ are equal to each other (i.e. the matrix $A$ does not change if we exchanges the places of these two rows) if and only if the vertices of $R_g$ corresponding to numbers $q$ and $r$ are equivalent according to relation $\sim$.

Analogous assertion is true about the columns of the matrix $A$ and the edges of the set $C_g$, which proves formula (\ref{bnk}).

From (\ref{mnbink1}) and (\ref{bnk}) it follows that
$$m!n! \sum_{g\in \mathfrak{G}_{m,n,k} } \frac{1}{\displaystyle \prod_{\delta \in [g]} \delta !} = \frac{(mn)!}{k!(mn-k)!} ,$$
which proves the theorem.

\hfill $\Box$

\bibliographystyle{plain}
\bibliography{umb2014}

\noindent Krasimir Yankov Yordzhev\\
South-West University ''N. Rilsky''\\
2700 Blagoevgrad, Bulgaria\\
Email: yordzhev@swu.bg

\end{document}